\documentclass [11pt] {article}

\usepackage{amsmath,amssymb,amsxtra,latexsym,amsthm}
\usepackage[english]{babel}
\usepackage{amsfonts}
\usepackage{multirow}
\usepackage{amssymb,amsmath}
\usepackage{amsfonts}
\usepackage{empheq,amsmath}
\usepackage{amsthm}
\usepackage{color}
\usepackage{colortbl}
\usepackage[table]{xcolor}
\usepackage{geometry}
\usepackage[numbers]{natbib}
\usepackage{mathrsfs}
\usepackage{stmaryrd}
\usepackage[normalem]{ulem}
\usepackage{enumerate}
\usepackage{graphicx}
\usepackage{float}
\usepackage{authblk}
\usepackage{ulem}

\begin {document}
\newtheorem{theorem}{Theorem}[section]
\newtheorem{proposition}[theorem]{Proposition}
\newtheorem{corollary}[theorem]{Corollary}
\newtheorem{lemma}[theorem]{Lemma}

\theoremstyle{definition}
\newtheorem{remark}[theorem]{Remark}

\thispagestyle{empty}

\title{On some generalized Fermat curves and chords of an affinely regular polygon inscribed in a hyperbola}

\author{
Herivelto Borges\thanks{{\it{Email address:}} hborges@icmc.usp.br}
\hspace{0.01cm} and
Mariana Coutinho\thanks{{\it{Email address:}} mariananery@alumni.usp.br}
}

\affil[]{Instituto de Ci\^encias Matem\'aticas e de Computa\c c\~ao\\ Universidade de S\~ao Paulo\\Avenida Trabalhador S\~ao-carlense, 400, CEP 13566-590, S\~ao Carlos SP, Brazil}

\date{}

\maketitle

\vspace{-1cm}

\begin{abstract}
Let $\mathcal{G}$ be the projective plane curve defined over $\mathbb{F}_q$ given by $$aX^nY^n-X^nZ^n-Y^nZ^n+bZ^{2n}=0,$$ where $ab\notin\{0,1\}$, and for each $s\in\{2,\ldots,n-1\}$, let $\mathcal{D}_s^{P_1,P_2}$ be the base-point-free linear series cut out on $\mathcal{G}$ by the linear system of all curves of degree $s$ passing through the singular points $P_1=(1:0:0)$ and $P_2=(0:1:0)$ of $\mathcal{G}$. The present work determines an upper bound for the number $N_q(\mathcal{G})$ of $\mathbb{F}_q$-rational points on the nonsingular model of $\mathcal{G}$ in cases where $\mathcal{D}_s^{P_1,P_2}$ is $\mathbb{F}_q$-Frobenius classical. As a consequence, when $\mathbb{F}_q$ is a prime field, the bound obtained for $N_q(\mathcal{G})$ improves in several cases the known bounds for the number $n_P$ of chords of an affinely regular polygon inscribed in a hyperbola passing through a given point $P$ distinct from its vertices.\\

\end{abstract}

%%%%%%%%%%%%%%%%%%%%%%%%%%%%%%%%%%%%%%%%%%%%%%%%%%%%%%%%%%
\section{Introduction}

Let $\mathbb{F}_q$ be the finite field with $q=p^m$ elements, where $p$ is a prime number, and let $\mathcal{X}$ be a (projective, nonsingular, geometrically irreducible, algebraic) curve of genus $g$ defined over $\mathbb{F}_q$. A fundamental problem in the theory of curves over finite fields is estimating the number $N_q(\mathcal{X})$ of $\mathbb{F}_q$-rational points on $\mathcal{X}$. Apart from a few classes of curves (see \cite{FCHN}, \cite{Moisio}), there is usually no explicit formula for $N_q(\mathcal{X})$. Nevertheless, some effective upper bounds for this number can be found in the literature. A famous example, given by the Hasse-Weil Theorem, is
\begin{eqnarray}
\label{HWbound}
N_q(\mathcal{X})\leqslant q+1+2gq^{1/2}.
\end{eqnarray}

Another noteworthy approach to bound $N_q(\mathcal{X})$ was established by St\"{o}hr and Voloch in $1986$ \cite{SV}. Their method, more geometric in nature, provides bounds that are dependent on the choice of an embedding of the curve in some $\mathbb{P}^M$, and which improve (\ref{HWbound}) in several circumstances (see \cite{GVFermat}, \cite{SV}).

For $a,b\in \mathbb{F}_q$ satisfying $ab\notin\{0,1\}$, let $\mathcal{G}$ be the projective plane curve with affine equation given by 
\begin{eqnarray}
\label{G}
g(X,Y)=aX^nY^n-X^n-Y^n+b=0.
\end{eqnarray}
This provides an example of a generalized Fermat curve, having recently been studied from the point of view of its automorphism group \cite{NP}.

The number $N_q(\mathcal{G})$ of $\mathbb{F}_q$-rational points on the nonsingular model $\mathcal{Y}$ of $\mathcal{G}$ was first investigated in the context of Finite Geometry to study the number of chords of an affinely regular polygon in $\mathbb{A}^2(\mathbb{F}_q)$ passing through a given point (see \cite{AK}, \cite{GIULIETTI}). 

A nondegenerate $k$-gon in the affine plane $\mathbb{A}^2(\mathbb{F}_q)$ is a set of $k$ pairwise distinct points arranged in a cyclic order in such a way that no three vertices are collinear. Here every $k$-gon is considered nondegenerate. If $A_1A_2\ldots A_k$ is a regular $k$-gon in the Euclidean plane, then a $k$-gon $B_1B_2\ldots B_k$ in the affine plane is affinely regular if the bijection $A_i\mapsto B_i$ preserves all parallelisms between chords (sides and diagonals). %Symbolically, $$A_{i_0}A_{j_0} \parallel A_{i_1}A_{j_1}\Leftrightarrow B_{i_0}B_{j_0} \parallel B_{i_1}B_{j_1},$$ for all $1\leqslant i_0<j_0\leqslant k$ and $1\leqslant i_1<j_1\leqslant k$. 

It is well-known that every affinely regular $k$-gon is inscribed in a conic, and for $p>2$, either $k\mid (q+1)$, or $k\mid (q-1)$, or $k=p$, according to whether the circumscribed conic is an ellipse, hyperbola, or parabola (see \cite{VCS}, \cite{K}, \cite{KS}). Moreover, if $k$ is large enough with respect to $q$, then the chords of any affinely regular $k$-gon cover most points of $\mathbb{A}^2(\mathbb{F}_q)$. The uncovered points are those remaining of the circumscribed conic and, in some cases, the center of the conic when this is either an ellipse or hyperbola (see \cite{K1}, \cite{KSS}, \cite{S}). This raises a natural question: what is the number $n_P$ of chords of an affinely regular $k$-gon of $\mathbb{A}^2(\mathbb{F}_q)$ passing through a given point $P$ distinct from its vertices?

In the particular case of such a $k$-gon being inscribed in an ellipse (resp. hyperbola), it was shown that determining $n_P$ is equivalent to determine the number of rational affine points (resp. the number of rational affine points lying in an appropriate subset of the plane) of a curve of the form of curve $\mathcal{G}$ (see \cite{AK}, \cite{GIULIETTI}). These relations give a connection between the problem of determining $n_P$ and that of studying $N_q(\mathcal{G})$. 

The primary method used in \cite{AK} and \cite{GIULIETTI} to give an upper bound for the number $N_q(\mathcal{G})$ was the St\"ohr-Voloch Theory. More precisely, in \cite{AK} Abatangelo and Korchm\'aros provided an upper bound for $N_q(\mathcal{G})$ based on the choice of an embedding of $\mathcal{Y}$ in  $\mathbb{P}^{5}$, and then on the study of its $\mathbb{F}_q$-Frobenius (non)classicality. Later on,  Giulietti  remarkably improved Abatangelo and Korchm\'aros' bound by considering a suitable embedding of $\mathcal{Y}$ in $\mathbb{P}^{3}$ \cite{GIULIETTI}. Further, more recently, and also using the St\"ohr-Voloch Theory, the number $N_q(\mathcal{G})$ was investigated in \cite{FG} in the context of generalized Fermat curves.

Accordingly, in terms of the St\"{o}hr-Voloch Theory, for $\overline{\mathbb{F}_q}(\mathcal{Y})=\overline{\mathbb{F}_q}(\mathcal{G})=\overline{\mathbb{F}_q}(x,y)$ and $s\in\{2,\ldots,n-1\}$, the present work determines an upper bound for $N_q(\mathcal{G})$ with respect to the morphism 
\begin{eqnarray}
\label{Morphism}
(\cdots:x^{i}y^{j}:\cdots):\mathcal{Y}\rightarrow \mathbb{P}^{N},
\end{eqnarray}
with $i,j$ integers satisfying $0\leqslant i,j \leqslant s-1$, $0\leqslant i+j\leqslant s$, and $N={{s+2}\choose{2}}-3$, in cases where it is $\mathbb{F}_q$-Frobenius classical. 

Further, based on techniques developed by Garcia and Voloch in [\citealp{GVFermat}, Section 3], and improved by Mattarei in \cite{Mattarei}, we focus our attention on the case where $q=p$ is prime. In particular, if $n=(p-1)/k\geqslant 3$ is a proper divisor of $p-1$, and $p \leqslant n^4/4$, for an affinely regular $k$-gon inscribed in a hyperbola, we obtain that $n_P$ is bounded roughly by 
\begin{eqnarray}
\label{bound}
3\cdot(2^{-1}\cdot k)^{2/3},
\end{eqnarray} 
which improves  in several cases the upper bound for $n_P$ given in [\citealp{GIULIETTI}, Theorem 4.1]:
\begin{eqnarray}
\label{BoundGiulietti}
n_P\leqslant \frac{k+1}{3}.
\end{eqnarray}

This paper is organized as follows. In Section \ref{Section2}, based on results of \cite{SV}, elements of the St\"ohr-Voloch Theory are recalled. In Section \ref{Section3}, the number $N_q(\mathcal{G})$ is studied. For each $s\in\{2,\ldots,n-1\}$, an upper bound for $N_q(\mathcal{G})$ is given  in cases where the morphism (\ref{Morphism}) is $\mathbb{F}_q$-Frobenius classical. Further, the case where $\mathbb{F}_q$ is the prime field $\mathbb{F}_p$ is addressed in Section \ref{Section3.1}. Finally, the bound given in (\ref{bound}) is presented in Section \ref{Section4}. 

%%%%%%%%%%%%%%%%%%%%%%%%%%%%%%%%%%%%%%%%%%%%%%%%%%%%%%%%%%
\begin{center}
\bf{Notation}
\end{center}

The following notation is used throughout this text.

\begin{itemize}
\item $\mathbb{F}_q$ is the finite field with $q=p^m$ elements, with $p$ a prime number.

\item $\overline{\mathbb{F}_q}$ is the algebraic closure of $\mathbb{F}_q$.

\item $P_1$ and $P_2$ are the points $(1:0:0)$ and $(0:1:0)$ of $\mathbb{P}^2(\overline{\mathbb{F}_q})$, respectively.

\item Unless otherwise stated, a curve denotes a projective, geometrically irreducible, algebraic curve.

\item For plane curves $\mathcal{F}$ and $\mathcal{F}'$, where $\mathcal{F}'$ is not necessarily irreducible and does not contain $\mathcal{F}$ as a component, $$\mathcal{F}'\cdot\mathcal{F}=\sum_{Q\in \mathcal{X}}I(Q,\mathcal{F}'\cap\eta)\,Q$$ is the intersection divisor cut out on $\mathcal{F}$ by $\mathcal{F}'$, where $\mathcal{X}$ is the nonsingular model of $\mathcal{F}$, and for each $Q\in\mathcal{X}$, $\eta$ is the corresponding branch of $\mathcal{F}$. Further, if $\mathcal{F}$ is defined over $\mathbb{F}_q$, then $N_q(\mathcal{F})$ is the number of $\mathbb{F}_q$-rational points on $\mathcal{X}$.

\item For a nonsingular curve $\mathcal{X}$ defined over $\mathbb{F}_q$, $\overline{\mathbb{F}_q}(\mathcal{X})$ is its function field, $\mathbb{F}_q(\mathcal{X})$ is its $\mathbb{F}_q$-rational function field, $\mathcal{X}(\mathbb{F}_q)$ is the set of its $\mathbb{F}_q$-rational points, and $N_q(\mathcal{X})$ is its number of $\mathbb{F}_q$-rational points.
\end{itemize}

%%%%%%%%%%%%%%%%%%%%%%%%%%%%%%%%%%%%%%%%%%%%%%%%%%%%%%%%%%
\section{Preliminaries}\label{Section2}

In this section, some elements of the St\"{o}hr-Voloch Theory based on \cite{SV} are recalled. 

Let $\mathcal{X}$ be a nonsingular curve of genus $g$ defined over $\mathbb{F}_q$. For a nondegenerate morphism $$\phi=(x_0:\cdots:x_M):\mathcal{X}\rightarrow \phi(\mathcal{X})\subset\mathbb{P}^M(\overline{\mathbb{F}_q}),$$ where $x_0,\ldots,x_M$ are functions in $\overline{\mathbb{F}_q}(\mathcal{X})$, let $\mathcal{D}$ be the corresponding base-point-free linear series of degree $\delta$ and dimension $M$ $$\bigg\{\text{div}\,(a_0x_0+\cdots+a_Mx_M)+E\,:\, (a_0:\cdots:a_M) \in \mathbb{P}^M(\overline{\mathbb{F}_q})\bigg\},$$ where $$E=\sum_{Q\in \mathcal{X}}e_QQ,$$ $$e_Q=-\min\{v_Q(x_0),\ldots,v_Q(x_M)\},$$ and $v_Q$ is the discrete valuation associated to the point $Q \in \mathcal{X}$. 

If $Q \in \mathcal{X}$, then $$\bigg\{v_Q(a_0x_0+\cdots+a_Mx_M)+e_Q\,:\, (a_0:\cdots:a_M) \in \mathbb{P}^M(\overline{\mathbb{F}_q})\bigg\}=\bigg\{j_0(Q),\ldots,j_M(Q)\bigg\},$$ where $0=j_0(Q)<\cdots<j_M(Q)$ and $v_Q(a_0x_0+\cdots+a_Mx_M)+e_Q=j_M(Q)$ for exactly one $(a_0:\cdots:a_M) \in \mathbb{P}^M(\overline{\mathbb{F}_q})$. Further, for almost all points $Q \in \mathcal{X}$ $$\{j_0(Q),\ldots,j_M(Q)\}=\{\epsilon_0,\ldots,\epsilon_M\},$$ where $0=\epsilon_0<\cdots<\epsilon_M$. The sequences $(j_0(Q),\ldots,j_M(Q))$ and $(\epsilon_0,\ldots,\epsilon_M)$ are called the $(\mathcal{D},Q)$-order sequence and the $\mathcal{D}$-order sequence, respectively, with $(\epsilon_0,\ldots,\epsilon_M)$ being also defined as the minimal sequence in the lexicographic order for which $$\det\,(D_t^{(\epsilon_i)}x_j)\neq 0,$$ where $t \in \overline{\mathbb{F}_q}(\mathcal{X})$ is a separable variable and $D_t^{(i)}$ is the $i$-th Hasse derivative with respect to $t$. Additionally, $\mathcal{D}$ (or $\phi$) is called classical if $(\epsilon_0,\ldots,\epsilon_M)=(0,\ldots,M)$, and nonclassical otherwise.

If $\phi$ is defined over $\mathbb{F}_q$, and $t \in \mathbb{F}_q(\mathcal{X})$, another important sequence related to $\mathcal{D}$ is the $\mathbb{F}_q$-Frobenius order sequence $$(\nu_0,\ldots,\nu_{M-1}),$$ which is the minimal sequence in the lexicographic order such that 
\begin{eqnarray*}
\label{definitionFrobeniusorders}
\displaystyle \det\left[\begin{array}{ccc} x_0^q&\cdots&x_M^q\\ D_t^{(\nu_0)}x_0&\cdots&D_t^{(\nu_0)}x_M\\\vdots&&\vdots \\D_t^{(\nu_{M-1})}x_0&\cdots&D_t^{(\nu_{M-1})}x_M\end{array}\right]\neq 0.
\end{eqnarray*}
This sequence satisfies $$\{\nu_0,\ldots,\nu_{M-1}\}=\{\epsilon_0,\ldots,\epsilon_M\}\setminus\{\epsilon_I\},$$ for some $I>0$. In this context, $\mathcal{D}$ (or $\phi$) is called $\mathbb{F}_q$-Frobenius classical if $(\nu_0,\ldots,\nu_{M-1})=(0,\ldots,M-1)$. Otherwise, it is called $\mathbb{F}_q$-Frobenius nonclassical.

The following result establishes an useful condition for the classicality and $\mathbb{F}_q$-Frobenius classicality of $\mathcal{D}$.

\begin{corollary}[\cite{SV}, Corollaries 1.8 and 2.7]
\label{condFqclass}
If $\delta<p$, then $\mathcal{D}$ is classical and $\mathbb{F}_q$-Frobenius classical.
\end{corollary}

In light of the previous considerations, this section ends with the following upper bound for the number $N_q(\mathcal{X})$, which is a refinement of [\citealp{SV}, Theorem 2.13] obtained from remarks at the beginning of [\citealp{SV}, Section 3].

\begin{theorem}[St\"ohr-Voloch]%----------------------------------------------------------------------------------------------------------------------
\label{SV}
\begin{eqnarray*}
N_q(\mathcal{X})\leqslant\frac{(\nu_0+\cdots+\nu_{M-1})\cdot(2g-2)+\delta\cdot(q+M)-\sum A(Q)}{M},
\end{eqnarray*}
where 
$$
A(Q)=\left\{\begin{array}{ll}
\displaystyle\sum_{i=1}^{M}(j_i(Q)-\nu_{i-1})-M,\mbox{ for {$Q \in \mathcal{X}(\mathbb{F}_q)$}}\\
\displaystyle\sum_{i=0}^{M-1}(j_i(Q)-\nu_i),\mbox{ otherwise.}
\end{array}\right.
$$
\end{theorem}%=====================================================================

%%%%%%%%%%%%%%%%%%%%%%%%%%%%%%%%%%%%%%%%%%%%%%%%%%%%%%%%%%
\section{The curve $\mathcal{G}$}\label{Section3}

For $a,b\in \mathbb{F}_q$ satisfying $ab\notin\{0,1\}$, let $\mathcal{G}$ be the plane curve defined over $\mathbb{F}_q$, with affine equation given by 
\begin{eqnarray*}
g(X,Y)=aX^nY^n-X^n-Y^n+b=0.
\end{eqnarray*}

The following result, which provides basic information about curve $\mathcal{G}$, can be found in [\citealp{AK}, Proposition 3.3].

\begin{proposition}%--------------------------------------------------------------------------------------------------------------------------------------
\label{curveG}
Let $n$ be a divisor of $q-1$. Then, the following holds:
\begin{enumerate}[\rm\indent1.]
\item $\mathcal{G}$ is geometrically irreducible.
\item The genus of $\mathcal{G}$ is $(n-1)^2$.
\item The only singular points of $\mathcal{G}$ are $P_1$ and $P_2$, and each singularity is ordinary with multiplicity $n$. Also, the tangent lines to $\mathcal{G}$ at $P_1$ and $P_2$ are given by the affine equations $Y=c$ and $X=c$, respectively, where $c^n=a^{-1}$, and those tangent lines intersect $\mathcal{G}$ at the corresponding points with multiplicity $2n$. 
\item The intersection multiplicity of a branch centered at $P_1$ or $P_2$ with its tangent line is $n+1$.
\item The points $(\xi:0:1)$ and $(0:\xi:1)$, with $\xi^n=b$, are inflection points of $\mathcal{G}$. Further, the tangent lines to $\mathcal{G}$ at $(\xi:0:1)$ and $(0:\xi:1)$ are given by the affine equations $X=\xi$ and $Y=\xi$, respectively, and those tangent lines intersect $\mathcal{G}$ at the corresponding points with multiplicity $n$.
\end{enumerate}
\end{proposition}%====================================================================

Let $\mathcal{Y}$ be the nonsingular model of $\mathcal{G}$ and $\overline{\mathbb{F}_q}(\mathcal{Y})=\overline{\mathbb{F}_q}(\mathcal{G})=\overline{\mathbb{F}_q}(x,y)$. For each $s\in\{2,\ldots,n-1\}$, consider the nondegenerate morphism $$\varphi_s^{P_1,P_2}=(\cdots:x^{i}y^{j}:\cdots):\mathcal{Y}\rightarrow \mathbb{P}^{N}(\overline{\mathbb{F}_q}),$$ with $i,j$ integers satisfying $0\leqslant i,j \leqslant s-1$ and $0\leqslant i+j\leqslant s$, which corresponds to the base-point-free linear series $\mathcal{D}_s^{P_1,P_2}\subset \text{Div}\,(\mathcal{Y})$ of dimension $N={{s+2}\choose{2}}-3$ cut out on $\mathcal{G}$ by all curves (not necessarily irreducible) of degree $s$ passing through the singularities $P_1$ and $P_2$ of $\mathcal{G}$.

The following proposition presents some important facts related to the linear series $\mathcal{D}_s^{P_1,P_2}$.

\begin{proposition}%--------------------------------------------------------------------------------------------------------------------------------------
\label{LinearSeries}
For each $s\in\{2,\ldots,n-1\}$, the following occurs:
\begin{enumerate}[\rm\indent1.]
\item \label{Item1} $\mathcal{D}_s^{P_1,P_2}$ has degree $\delta=2n\cdot(s-1)$.
\item \label{Item2} For $Q\in \mathcal{Y}$ corresponding to a point of $\mathcal{G}$ of the form $(\xi:0:1)$ or $(0:\xi:1)$, with $\xi^n=b$, the $(\mathcal{D}_s^{P_1,P_2},Q)$-order sequence is given by the elements of $$\bigg\{i+jn\,:\,0\leqslant i,j\leqslant s-1 \text{ and } 0\leqslant i+j\leqslant s\bigg\}.$$
\item \label{Item3} For $Q_\eta\in \mathcal{Y}$ corresponding to a branch $\eta$ of $\mathcal{G}$ centered at $P_1$ or $P_2$, the $(\mathcal{D}_s^{P_1,P_2},Q_\eta)$-order sequence is given by the elements of $$\displaystyle\bigg\{i+j(n+1)-1\,:\,0\leqslant i,j \text{ and } 0\leqslant i+j\leqslant s\bigg\}\setminus\bigg\{-1,s(n+1)-1\bigg\}.$$ 
\end{enumerate}
\end{proposition}
\begin{proof}
For each $s\in\{2,\ldots,n-1\}$, let $\Sigma^{P_1,P_2}_s$ be the linear system of all curves (not necessarily irreducible) of degree $s$ passing through the singularities $P_1$ and $P_2$ of $\mathcal{G}$, and let $$\tilde{\mathcal{D}_s}^{P_1,P_2}=\bigg\{\mathcal{F}\cdot\mathcal{G}\,:\,\mathcal{F}\in\Sigma^{P_1,P_2}_s \bigg\}$$ be the linear series cut out on $\mathcal{G}$ by the linear system $\Sigma^{P_1,P_2}_s$. By B\'ezout's Theorem, $\tilde{\mathcal{D}_s}^{P_1,P_2}$ has degree $2ns$.

The base locus of $\tilde{\mathcal{D}_s}^{P_1,P_2}$ is the divisor $$\ell_{\infty}\cdot \mathcal{G}=Q_1^{(1)}+\cdots+Q_n^{(1)}+Q_1^{(2)}+\cdots +Q_n^{(2)}$$ of degree $2n$, where the $Q_j^{(i)}$'s are all the distinct points on $\mathcal{Y}$ for which the corresponding branches are centered at $P_i$, for $i=1,2$ and $j=1,\ldots,n$, and $\ell_{\infty}$ is the line given by the equation $Z=0$. 

Therefore, 
\begin{eqnarray}
\label{DefLinSer}
\mathcal{D}_s^{P_1,P_2}=\tilde{\mathcal{D}_s}^{P_1,P_2}-\ell_{\infty}\cdot \mathcal{G}=\bigg\{D-\ell_{\infty}\cdot \mathcal{G}\,:\,D\in\tilde{\mathcal{D}_s}^{P_1,P_2}\bigg\}
\end{eqnarray}
 has degree $2ns-2n=2n\cdot(s-1)$, which proves statement \ref{Item1}.

To prove statement \ref{Item2}, let $Q \in \mathcal{Y}$ be a point corresponding to  $Q'\in\mathcal{G}$, where $Q'$ is equal to $(\xi:0:1)$ or $(0:\xi:1)$, with $\xi^n=b$, and let $\ell_i=\overline{Q'P_i}$, for $i=1,2$. One can verify that none of the lines $\ell_1$ and $\ell_2$ is equal to $\ell_{\infty}$, and exactly one of them is the tangent line to $\mathcal{G}$ at $Q'$. Further, from Proposition \ref{curveG}, the $(\mathcal{D}_1,Q)$-order sequence is $(0,1,n)$, where $\mathcal{D}_1$  is the linear series cut out on $\mathcal{G}$ by the linear system of lines. Therefore, considering the reducible curves $\mathcal{C}$ given by the union of $s$ lines chosen (with multiplicity) in the set $\{\ell_1,\ell_2,\ell_{\infty}\}$, the $(\mathcal{D}_s^{P_1,P_2},Q)$-order sequence is given by the elements of $$\bigg\{i+jn\,:\,0\leqslant i,j\leqslant s-1 \text{ and } 0\leqslant i+j\leqslant s\bigg\}.$$ 

Finally, if $Q_\eta\in \mathcal{Y}$ corresponds to a branch $\eta$ of $\mathcal{G}$ centered at $P_1$ or $P_2$, from Proposition \ref{curveG}, the $(\mathcal{D}_1,Q_\eta)$-order sequence is $(0,1,n+1)$. Thus, considering the reducible curves $\mathcal{C}$ passing through $P_1$ and $P_2$, and given by the union of $s$ lines (possibly chosen with multiplicity), equation (\ref{DefLinSer}) shows that the $(\mathcal{D}_s^{P_1,P_2},Q_\eta)$-order sequence is given by the elements of $$\bigg\{i+j(n+1)-1\,:\,0\leqslant i,j \text{ and } 0\leqslant i+j\leqslant s\bigg\}\setminus\bigg\{-1,s(n+1)-1\bigg\},$$ which completes the proof.
\end{proof}%=======================================================================

Based on Theorem \ref{SV} and Proposition \ref{LinearSeries}, the following result gives an upper bound for the number $N_q(\mathcal{G})$ in cases where $\mathcal{D}_s^{P_1,P_2}$ is $\mathbb{F}_q$-Frobenius classical.

\begin{corollary}%-----------------------------------------------------------------------------------------------------------------------------------------
\label{Corollary}
Let $s\in\{2,\ldots,n-1\}$. If $\mathcal{D}_s^{P_1,P_2}$ is $\mathbb{F}_q$-Frobenius classical, then 
\begin{eqnarray}
\label{boundSVwR}
N_q(\mathcal{G})&\leqslant&(N-1)\cdot(n^2-2n)+\frac{\delta\cdot(q+N)}{N}-2\cdot\frac{n_1\cdot \alpha +n_2\cdot \beta+n\cdot \gamma}{N},
\end{eqnarray}
 where
\begin{itemize}
\item $n_1$ and $n_2$ are the number of roots in $\mathbb{F}_q$ of the polynomials $T^n-b$ and $T^n-a^{-1}$, respectively
\item $\alpha=1+(s-1)\cdot n-N$
\item $\beta=(s-1)\cdot(n+1)-N $
\item $\gamma=2(n+1)-s\cdot(4n+3)-N\cdot(N-1)+(s\cdot(2n+3)-3)\cdot\frac{N+3}{3}$.
\end{itemize}
\end{corollary}

\begin{proof}
The part $$(N-1)\cdot(n^2-2n)+\frac{\delta\cdot(q+N)}{N}$$ of (\ref{boundSVwR}) follows directly from Theorem \ref{SV}, since the genus of $\mathcal{G}$ is $(n-1)^2$ by Proposition \ref{curveG}.

Further, let $Q\in \mathcal{Y}$ be a point corresponding to $(\xi:0:1)$ or $(0:\xi:1) \in \mathcal{G}$, with $\xi^n=b$, and let $Q_\eta\in\mathcal{Y}$ be a point corresponding to a branch $\eta$ of $\mathcal{G}$ centered at $P_1$ or $P_2$. Using the notation as in Theorem \ref{SV}, from Proposition \ref{LinearSeries}, the numbers $A(Q)$ and $A(Q_\eta)$ are given by the following expressions:
\begin{empheq}[left={A(Q)=\empheqlbrace}]{align}
&s\cdot (n+1)\cdot \frac{-6+(s+1)\cdot(s+2)}{6}-\frac{N\cdot(N-1)}{2}-N,\text{ if } Q\in \mathcal{Y}(\mathbb{F}_q)\label{AP1}\\
&s\cdot (n+1)\cdot \frac{-6+(s+1)\cdot(s+2)}{6}-\frac{N\cdot(N-1)}{2}-(1+(s-1)\cdot n),\text{ otherwise}\label{AP2}
\end{empheq}

\noindent and 
\begin{empheq}[left={A(Q_\eta)=\empheqlbrace}]{align}
&2-s\cdot(n+1)+(s\cdot(n+2)-3)\cdot\frac{(s+1)\cdot (s+2)}{6}-\frac{N(N-1)}{2}-N,\text{ if } Q_\eta\in \mathcal{Y}(\mathbb{F}_q)\label{APgamma1}\\
&2-s\cdot(n+1)+(s\cdot(n+2)-3)\cdot\frac{(s+1)\cdot (s+2)}{6}-\frac{N(N-1)}{2}-(s-1)\cdot(n+1),\text{ otherwise}\label{APgamma2}
\end{empheq}

\noindent since 
\begin{eqnarray*}
\sum_{i=1}^{N} j_i(Q)=s\cdot (n+1)\cdot \frac{-6+(s+1)\cdot(s+2)}{6}\, \text{ and }\,j_N(Q)=1+(s-1)\cdot n,
\end{eqnarray*}

\noindent and
\begin{eqnarray*}
\sum_{i=1}^{N} j_i(Q_\eta)=2-s\cdot(n+1)+(s\cdot(n+2)-3)\cdot\frac{(s+1)\cdot (s+2)}{6}\,\text{ and }\,j_N(Q_\eta)=(s-1)\cdot(n+1).                                          
\end{eqnarray*}

\noindent Therefore,
\begin{eqnarray*}
N_q(\mathcal{G})&\leqslant&(N-1)\cdot(n^2-2n)+\frac{\delta\cdot(q+N)}{N}\\
&&-\frac{2n_1\cdot A(Q)^{(\ref{AP1})}+2(n-n_1)\cdot A(Q)^{(\ref{AP2})}+2n_2\cdot A(Q_\eta)^{(\ref{APgamma1})}+2(n-n_2)\cdot A(Q_\eta)^{(\ref{APgamma2})}}{N}\\
&=&(N-1)\cdot(n^2-2n)+\frac{\delta\cdot(q+N)}{N}-2\cdot\frac{n_1\cdot \alpha +n_2\cdot \beta+n\cdot \gamma}{N},
\end{eqnarray*}
 where\\
 
\indent $\bullet\,$ $n_1$ and $n_2$ are the number of roots in $\mathbb{F}_q$ of the polynomials $T^n-b$ and $T^n-a^{-1}$, respectively\\

\indent $\bullet\,$ $\alpha=A(Q)^{(\ref{AP1})}-A(Q)^{(\ref{AP2})}=1+(s-1)\cdot n-N$\\

\indent $\bullet\,$ $\beta=A(Q_\eta)^{(\ref{APgamma1})}-A(Q_\eta)^{(\ref{APgamma2})}=(s-1)\cdot(n+1)-N $\\

\indent $\bullet\,$ $\gamma=A(Q)^{(\ref{AP2})}+A(Q_\eta)^{(\ref{APgamma2})}=2(n+1)-s\cdot(4n+3)-N\cdot(N-1)+(s\cdot(2n+3)-3)\cdot\frac{N+3}{3}$.
\end{proof}%=======================================================================

%%%%%%%%%%%%%%%%%%%%%%%%%%%%%%%%%%%%%%%%%%%%%%%%%%%%%%%%%%
\subsection{The case $q=p$}\label{Section3.1}

Hereafter, let $\mathbb{F}_q$ be the prime field $\mathbb{F}_p$. For each $u\in [2,+\infty)$ and $t_0\in[6,+\infty)$, let us consider
$$
\begin{array}{cccc}
f_u:&[6,+\infty)&\rightarrow&\mathbb{R}\\
&t&\mapsto&\frac{3t^2-23t+26}{6}+ 4\cdot\frac{u + 3}{t}
\end{array}, 
$$

\noindent and 
\begin{eqnarray}
k_{t_0}:=\frac{t_0\cdot(t_0+1)\cdot(3t_0-10)}{12}-3=\frac{1}{4}t_0^3 - \frac{7}{12}t_0^2 - \frac{5}{6}t_0 - 3.
\end{eqnarray}

An important step in the proof of Theorem \ref{Bound} is the following lemma, whose proof is straightforward.

\begin{lemma}%--------------------------------------------------------------------------------------------------------------------------------------------
\label{Lema1}
Let $t_0\in [6,+\infty)$. For $u\in[2,+\infty)$,
\begin{eqnarray}
\label{fundequal}
u\leqslant k_{t_0}&\text{if and only if}&f_u(t_0)\leqslant f_u(t_0+1).	
\end{eqnarray}
\end{lemma}%=======================================================================

\begin{theorem}%--------------------------------------------------------------------------------------------------------------------------------------
\label{Bound}
Let $n\geqslant 3$ be a proper divisor of $p-1$ and $k=(p-1)/n$. If $p-1\leqslant n\cdot\bigg(\frac{(n+3)\cdot(n+4)\cdot(3n-1)}{12}-3\bigg)$, then
\begin{eqnarray}
\label{BoundNq}
N_p(\mathcal{G})\leqslant n^2\cdot\bigg(3\cdot(2^{1/2}\cdot k)^{2/3}-\frac{103}{19}\cdot(2^{1/2}\cdot k)^{1/3}+\frac{13}{3}\bigg).
\end{eqnarray}
\end{theorem}

\begin{proof}
From Corollary \ref{Corollary} and its proof,
\begin{eqnarray*}
N_p(\mathcal{G})&\leqslant&(N-1)\cdot(n^2-2n)+\frac{\delta\cdot(p+N)}{N}-2\cdot \frac{n_1\cdot \alpha +n_2\cdot \beta+n\cdot \gamma}{N}\\
&\leqslant&(N-1)\cdot(n^2-2n)+\frac{\delta\cdot(p+N)}{N}-\frac{2n\cdot \gamma}{N}\\
&\leqslant&\frac{\frac{1}{2}s^3 + \frac{13}{6}s^2 - \frac{7}{3}s - 4}{s + 4}\cdot n^2 +\frac{4n\cdot((p-1)\cdot(s-1)+n)}{(s + 4)\cdot(s-1)}\\
&\leqslant&n^2\cdot\bigg(\frac{\frac{1}{2}s^3 + \frac{13}{6}s^2 - \frac{7}{3}s - 4}{s + 4}+4\cdot\bigg(\frac{k+1}{s+4}\bigg)\bigg),
\end{eqnarray*}
for each $s\in\{2,\ldots,n-1\}$ satisfying $2n\cdot(s-1)\leqslant p-1$, since the latter condition, together with Corollary \ref{condFqclass} and Proposition \ref{LinearSeries},  implies that $\mathcal{D}_s^{P_1,P_2}$ is $\mathbb{F}_{p}$-Frobenius classical. Hence, defining $t:=s+4$, it follows that
$$N_p(\mathcal{G})\leqslant n^2\cdot\bigg(\frac{3t^2-23t+26}{6}+ 4\cdot\frac{k + 3}{t}\bigg),$$
for all $t\in\{6,\ldots,n+3\}$ satisfying $t\leqslant \frac{k}{2}+5$, and
$$\frac{N_p(\mathcal{G})}{n^2}\leqslant V(k):= \min\bigg\{\frac{3t^2-23t+26}{6}+ 4\cdot\frac{k + 3}{t}\,:\, t\in\{6,\ldots,n+3\} \text{ and } t\leqslant \frac{k}{2}+5\bigg\}.$$

One can check that $k<k_{t_0}$, for all integers $t_0\geqslant \frac{k}{2}+5$. Thus, an iterated application of Lemma \ref{Lema1} yields $$V(k)=\min\{f_k(t)\,:\, t\in\mathbb{Z}\cap[6,n+3]\}.$$ 

Likewise, from the assumption $k\leqslant k_{n+3}=\frac{(n+3)\cdot(n+4)\cdot (3n-1)}{12}-3$, we further have $$V(k)=\min\{f_{k}(t)\,:\, t\in \mathbb{Z}\cap[6,+\infty)\}.$$

The objective is to determine a suitable function $W$ such that $W(k)\geqslant V(k)$, for all prime numbers $p$ and for all proper divisors $n\geqslant 3$ of $p-1$, where $k=(p-1)/n$. Then, defining for each $u \in[2,+\infty)$ $$\tilde{V}(u):=\min\{f_{u}(t)\,:\, t\in \mathbb{Z}\cap[6,+\infty)\},$$ Lemma \ref{Lema1} implies that 
$$
\tilde{V}(u)=\left\{
\begin{array}{ll}
f_u(t_0+1),&\text{if } k_{t_0}\leqslant u\leqslant k_{t_0+1}\text{ for some }t_0\in\mathbb{Z}\cap[6,+\infty)\\
f_u(6),&\text{if } 2\leqslant u\leqslant k_6
\end{array}\right.,
$$
and that $$\tilde{V}(k_{t_0})=f_{k_{t_0}}(t_0)=\frac{3}{2}t_0^2 - \frac{37}{6}t_0 + 1.$$ 

Thus one may choose a convenient concave function $W$ such that $W(k_{t_0})\geqslant \tilde{V}(k_{t_0})$ for all $t_0\in\mathbb{Z}\cap[6,+\infty)$, and such that $W(u)\geqslant \tilde{V}(u)$, for all integers  $2\leqslant u< k_6$. Indeed, fixed $u\in[k_6,+\infty)$, if $k_{t_0}\leqslant u\leqslant k_{t_0+1}$, for some $t_0\in\mathbb{Z}\cap[6,+\infty)$, let $u=u_1k_{t_0}+u_2k_{t_0+1}$, where $0\leqslant u_1,u_2\leqslant 1$ and $u_1+u_2=1$. Then
\begin{eqnarray*}
W(u)&=&W(u_1k_{t_0}+u_2k_{t_0+1})\\
&\geqslant&u_1W(k_{t_0})+u_2W(k_{t_0+1})\text{ (by the concavity of function }W)\\
&\geqslant&u_1\tilde{V}(k_{t_0})+u_2\tilde{V}(k_{t_0+1})\\
&=&\tilde{V}(u).
\end{eqnarray*}

Let us consider the family of concave functions $W_{\lambda}$ defined in $[2,+\infty)$ by $$W_{\lambda}(u)=3\cdot(2^{1/2}\cdot u)^{2/3}-\frac{103}{19}\cdot(2^{1/2}\cdot u)^{1/3}+\lambda,$$ where $\lambda$ is a constant. 

Direct computation shows that the smallest $\lambda$ such that $W_{\lambda}(k_{t_0})\geqslant \tilde{V}(k_{t_0})$ for all $t_0\in\mathbb{Z}\cap[6,+\infty)$, and such that $W_{\lambda}(u)\geqslant \tilde{V}(u)=f_u(6)$ for all integers  $2\leqslant u<k_6$, is $$\frac{10}{3}-\frac{103}{19}\cdot 2^{1/2}\cong \frac{13}{3}.$$ 

Therefore, the proof completes considering $W$ in $[2,+\infty)$ defined as $W(u)=3\cdot(2^{1/2}\cdot u)^{2/3}-\frac{103}{19}\cdot(2^{1/2}\cdot u)^{1/3}+\frac{13}{3}.$
\end{proof}%=======================================================================

%%%%%%%%%%%%%%%%%%%%%%%%%%%%%%%%%%%%%%%%%%%%%%%%%%%%%%%%%%
\section{Number of chords of an affinely regular polygon inscribed in a hyperbola passing through a given point}\label{Section4}

Let $A=A_1A_2\ldots A_k$ be an affinely regular polygon with $k=(p-1)/n$ vertices inscribed in a hyperbola $\mathcal{H}$ of $\mathbb{A}^2(\mathbb{F}_p)$. Changing variable, one may consider $\mathcal{H}$ given by the equation $XY=1$.

For $a,b\in \mathbb{F}_p$ satisfying $ab\notin\{0,1\}$, let $\mathcal{G}$ be the plane curve defined over $\mathbb{F}_p$, with affine equation
\begin{eqnarray*}
g(X,Y)=aX^nY^n-X^n-Y^n+b=0.
\end{eqnarray*}

The following result  provides a relation between the number of chords of $A$ passing through the point $P=(a,b)\in\mathbb{A}^2(\mathbb{F}_p)$ and the number $N_p(\mathcal{G})$.

\begin{proposition}[\cite{GIULIETTI}, Proposition 2.1]%-----------------------------------------------------------------------------------------
\label{RelationNn}
Let $N_p$ be the number of $\mathbb{F}_p$-rational affine points of $\mathcal{G}$ not lying in the coordinate axes or the line of equation $X=Y$. Then $n_P=\frac{N_p}{2n^2}$.
\end{proposition}%====================================================================

Therefore, from Theorem \ref{Bound} and Proposition \ref{RelationNn}, the following upper bound for $n_P$ is obtained.

\begin{corollary}%-----------------------------------------------------------------------------------------------------------------------------------------
\label{Corn_P}
If $n\geqslant 3$ is a proper divisor of $p-1$ such that $p-1\leqslant n\cdot\bigg(\frac{(n+3)\cdot(n+4)\cdot(3n-1)}{12}-3\bigg)$, then
\begin{eqnarray}
n_P\leqslant \frac{1}{2}\cdot\bigg(3\cdot(2^{1/2}\cdot k)^{2/3}-\frac{103}{19}\cdot(2^{1/2}\cdot k)^{1/3}+\frac{13}{3}\bigg).\label{Boundnp}
\end{eqnarray}
\end{corollary}%=====================================================================

\begin{remark}
One can check that bound ({\ref{Boundnp}}) improves that given in ({\ref{BoundGiulietti}}) for $k\geqslant 44$. On the other hand, for each $2\leqslant k<44$, bounds (\ref{BoundGiulietti}) and (\ref{Boundnp}) effectively differ by at most one unit. Moreover, considering the definition of $V(k)=\tilde{V}(k)$ as given in the proof of Theorem {\ref{Bound}}, for $25<k<44$, we obtain that $n_P\leqslant \lfloor \frac{1}{2}V(k)\rfloor \leqslant \lfloor\frac{k+1}{3}\rfloor$, with $\frac{k+1}{3}$ as in ({\ref{BoundGiulietti}}). For an illustration of the region where ({\ref{Boundnp}}) is better than the bound for $n_P$ derived from the Hasse-Weil bound, see also Figure {\ref{Graph}}.

\begin{figure}[H]
\begin{center}
\includegraphics[scale=0.7]{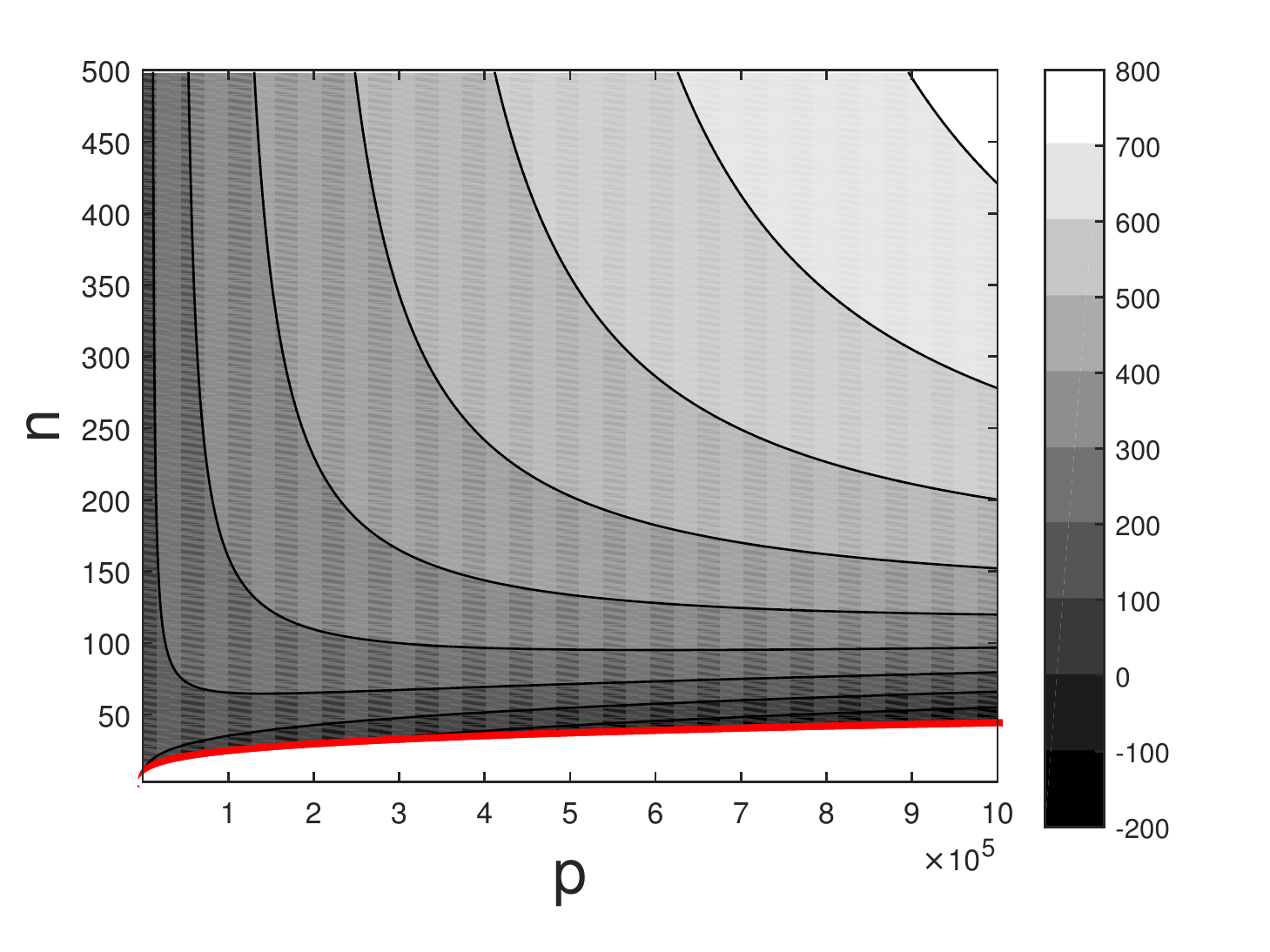}
\caption{Contour plot of the difference $\Delta(p,n)=(\frac{p + 1 + 2\cdot (n - 1)^2\cdot p^{1/2}}{2\cdot n^2})- \frac{1}{2}\cdot(3\cdot (2^{1/2}\cdot k)^{2/3} - \frac{103}{19}\cdot(2^{1/2}\cdot k)^{1/3} + \frac{13}{3})$ between the bounds for the number $n_P$ obtained from the Hasse-Weil bound and Corollary \ref{Corn_P}, for $n< p-1\leqslant n \cdot (\frac{(n+3)\cdot(n+4)\cdot (3n-1)}{12}-3)$ and $k=(p-1)/n$. In red, the curve $p-1=n \cdot (\frac{(n+3)\cdot(n+4)\cdot (3n-1)}{12}-3)$.\label{Graph}}
\end{center}
\end{figure}
\end{remark}

%%%%%%%%%%%%%%%%%%%%%%%%%%%%%%%%%%%%%%%%%%%%%%%%%%%%%%%%%%
\section*{Aknowlegments}

The first author was supported by FAPESP (Brazil), grant 2017/04681-3. The study of the second author was financed in part by the Coordena\c{c}\~{a}o de Aperfei\c{c}oamento de Pessoal de N\'{i}vel Superior - Brasil (CAPES) - Finance Code 001, and CNPq (Brazil), grant 154359/2016-5.

%%%%%%%%%%%%%%%%%%%%%%%%%%%%%%%%%%%%%%%%%%%%%%%%%%%%%%%%%%


\begin{thebibliography}{99}
\bibitem{AK} ABATANGELO, V.; KORCHM\'{A}ROS, G. Una generalizzazione di un teorema di B. Segre sui punti regolari rispetto ad una ellisse di un piano affine di Galois. {\it{Annali di Matematica Pura ed Applicata}}, v. 172, p. 87--102, 1997.
%\bibitem{NRBSPC} ARAKELIAN, N. Number of rational branches of a singular plane curve over a finite field. {\it{Finite Fields and their Applications}}, v. 48, p. 87--102, 2017.
%\bibitem{CFCCFG} ARAKELIAN, N.; BORGES, H. Frobenius nonclassicality with respect to linear systems of curves of arbitrary degree. {\it{Acta Arithmetica}}, v. 167, p. 43--66, 2015.
\bibitem{FCHN} ARAKELIAN, N.; BORGES, H. Frobenius nonclassicality of Fermat curves with respect to cubics. {\it{Israel Journal of Mathematics}}, v. 218, p. 273--297, 2017.
\bibitem{NP} ARAKELIAN, N.; SPEZIALI, P. On generalizations of Fermat curves over finite fields and their automorphisms. {\it{Communications in Algebra}}, v. 45, p. 4926--4938, 2017. 
%\bibitem{CV} CARLIN, M. L.; VOLOCH, J. F. Plane curves with many points over finite fields. {\it{Rocky Mountain Journal of Mathematics}}, v. 34, p. 1255--1259, 2004.
%\bibitem{COULTER} COULTER, R. S. Explicit evaluations of some Weil sums. {\it{Acta Arithmetica}}, v. 83, p. 241--251, 1998.
\bibitem{VCS} VAN DE CRAATS, J.; SIMONIS, J. Affinely regular polygons. {\it{ Nieuw Archief voor Wiskunde}}, v. IV, p. 225--240, 1986.
\bibitem{FG} FANALI, S.; GIULIETTI, M. On the number of rational points of generalized Fermat curves over finite fields. {\it{International Journal of Number Theory}}, v. 8, p. 1087--1097, 2012.
%\bibitem{FMac} FRIED, M. D.; MacRAE, R. E. On curves with separated variables. {\it{Math. Ann.}}, v. 180, p. 220--226, 1969.
%\bibitem{FULTON} FULTON, W. {\textit{Algebraic Curves: An Introduction to Algebraic Geometry}}. 3. ed. 2008.
%\bibitem{AGYL} GARCIA, A. L. P.; LEQUAIN, Y. A. E. {\textit{Elementos de \'{a}lgebra}}. 6. ed. Rio de Janeiro: IMPA, 2012.
%\bibitem{GVWronskians} GARCIA, A.; VOLOCH, J. F. Wronskians  and linear independence in fields of prime characteristic. {\it{Manuscripta Math.}}, v. 59, p. 457--469, 1987.
\bibitem{GVFermat} GARCIA, A.; VOLOCH, J. F. Fermat curves over finite fields. {\it{Journal of Number Theory}}, v. 30, p. 345--356, 1988.
\bibitem{GIULIETTI} GIULIETTI, M. On the number of chords of an affinely regular polygon passing through a given point. {\it{Acta Scientiarum Mathematicarum (Szeged)}}, v. 74, p. 901--913, 2008.
%\bibitem{GOPPA} GOPPA, V. D. {\textit{Geometry and codes}}. 1 ed. Dordrecht: Kluwer Academic Publishers, 1988.
%\bibitem{H} HIRSCHFELD, J. W. P. {\textit{Projective geometries over finite fields}}. 1. ed. Oxford: Oxford University Press, 1979.
%\bibitem{HK} HIRSCHFELD, J. W. P.; KORCHM\'{A}ROS, G. On the number of solutions of an equation over a finite field. {\textit{Bull. London Math. Soc.}}, v. 33, p. 16--24, 2001.
%\bibitem{HKT} HIRSCHFELD, J. W. P.; KORCHM\'{A}ROS, G.; TORRES, F. {\textit{Algebraic Curves over a Finite Field}}. 1. ed. Princeton: Princeton University Press, 2008.
\bibitem{K} KORCHM\'{A}ROS, G. Poligoni affini regolari dei piani di Galois di ordine dispari. {\it{Atti della Accademia Nazionale dei Lincei. Classe di Scienze Fisiche, Matematiche e Naturali. Rendiconti Lincei}}, v. 56, p. 690--697, 1974.
\bibitem{K1} KORCHM\'{A}ROS, G. New examples of $k$-arcs in $PG(2, q)$. {\it{European Journal of Combinatorics}}, v. 4, p. 329--334, 1983.
\bibitem{KSS} KORCHM\'{A}ROS, G.; STORME, L.; SZ\H{O}NYI, T. Space-filling subsets of a normal rational curve. {\it{Journal of Statistical Planning and Inference}}, v. 58, p. 93--110, 1997.
\bibitem{KS} KORCHM\'{A}ROS, G.; SZ\H{O}NYI, T. Affinely regular polygons in an affine plane. {\it{Contributions to Discrete Mathematics}}, v. 3, p. 20--38, 2008.
\bibitem{Mattarei} MATTAREI, S. On a bound of Garcia and Voloch for the number of points of a Fermat curve over a prime field. {\it{Finite Fields and their Applications}}, v. 13, p. 773--777, 2007.
%\bibitem{McEliece} McELIECE, R. J. {\textit{Finite Fields for Computer Scientists and Engineers}}. Dordrecht: Kluwer, 1987.
\bibitem{Moisio} MOISIO, M. On the number of rational points on some families of Fermat curves over finite fields. {\it{Finite Fields and their Applications}}, v. 13, p. 546--562, 2007.
%\bibitem{RP} PARDINI, R. Some remarks on plane curves over fields of finite characteristic. {\it{Compositio Math.}}, v. 60, p. 3--17, 1986.
%\bibitem{HS} STICHTENOTH, H. Algebraic Function Fields and Codes. 2. ed. New York: Springer-Verlag, 2009.
\bibitem{SV} ST\"OHR, K. O.; VOLOCH, J. F. Weierstrass points and curves over finite fields. {\textit{Proceedings of the London Mathematical Society}}, v. 52, n. 1, p. 1--19, 1986.
\bibitem{S} SZ\H{O}NYI, T. Note on the order of magnitude of $k$ for complete $k$-arcs in $PG(2, q)$. {\it{Discrete Mathematics}}, v. 66, p. 279--282, 1987.
%\bibitem{FT} TORRES, F. The Approach of St\"{o}hr-Voloch to the Hasse-Weil Bound with Applications to Optimal Curves and Plane Arcs. {\it{https://arxiv.org/abs/math/0011091v1}}
\end{thebibliography}
\end {document}